\documentclass[12pt,a4paper,reqno]{article}
\usepackage[scaled]{helvet}

\usepackage{graphicx}
\usepackage{epstopdf} 
\usepackage{amsfonts}
\usepackage{amsthm}
\usepackage{amsmath}
\usepackage{amscd}
\usepackage[mathscr]{eucal}
\usepackage{indentfirst}

\usepackage{pict2e}
\usepackage{epic}
\numberwithin{equation}{section}
\usepackage{MnSymbol}
\usepackage[top=1.5cm, bottom=1.5cm, left=1.5 cm, right=1cm]{geometry}
\usepackage{tikz}
\usepackage{enumitem}
\usepackage{bm}

\usepackage{hyperref}
\numberwithin{equation}{section}

\theoremstyle{plain}
\newtheorem{Th}{Theorem}[section]
\newtheorem{Lemma}[Th]{Lemma}
\newtheorem{Cor}[Th]{Corollary}

 \theoremstyle{definition}
\newtheorem{Def}[Th]{Definition}

\newtheorem{?}[Th]{Problem}
\newtheorem{Ex}[Th]{Example}
\newcommand{\R}{\mathbb{R}}
\newcommand\Q{\mathbb{H}}
\newcommand\F{\mathcal{F}}
\newcommand \be     {\begin{equation}}
\newcommand \ee     {\end{equation}}
\newcommand \bee     {\begin{eqnarray*}}
\newcommand \eee     {\end{eqnarray*}}
\newcommand{\x}{\bm{x}}
\newcommand{\y}{\bm{y}}
\newcommand{\xii}{\bm{\xi}}
\newcommand{\e}{\bm{i}}
\newcommand{\f}{\bm{j}}
\newcommand{\g}{\bm{k}}
\newcommand\myeq{\mathrel{\overset{\makebox[0pt]{\mbox{\normalfont\tiny\sffamily def}}}{=}}}

\begin{document}
\newpage
\vskip .2cm 
\begin{center}
\huge\textbf{Miyachi's Theorem For the Quaternion Fourier Transform}
\end{center}

\begin{center}
\textbf{Youssef El Haoui$^{1,}$\footnote[1]{Corresponding author.}, Said Fahlaoui$^1$}

 $^1$Department of Mathematics and Computer Sciences, Faculty of Sciences, Equipe d’Analyse Harmonique et Probabilités, University Moulay Ismail, BP 1120, Zitoune, Meknes, Morocco\\

E-MAIL: y.elhaoui@edu.umi.ac.ma, s.fahlaoui@fs.umi.ac.ma
\end{center}

\vspace{1cm}
\begin{abstract}
	The quaternion Fourier transform (QFT) satisfies some uncertainty principles similar to the Euclidean Fourier transform. In this paper, we establish Miyachi's theorem for this transform and consequently generalize and prove the analogue of Hardy's theorem and Cowling-Price uncertainty principle in the QFT domain. 
\end{abstract}



	\section{Introduction} 
	Uncertainty principle (UP) is an important tool in harmonic anlysis; it states that a nonzero function and its Fourier transform cannot both be very rapidly decreased. UP has implications in different areas like quantum physics, information processing, signal analysis, etc.\\
	In signal analysis, it gives in general a lower bound for the simultaneous localization of signals in phase and frequency spaces.
There are many advantageous ways to get the statement about localization precise; examples include theorems of Hardy theorem \cite{HA33}, Beurling \cite{HO91}, and Miyachi \cite{MI97} which interpreted differently the localisation, as sharp pointwise estimates of a signal and its Fourier transform. More precisely Miyachi's theorem asserts that if $f$ is a measurable function on $\R$ satisfying :
	
\[\begin{aligned} e^{{\alpha x}^2}f\in L^1({{\mathbb {R}}})+ L^{\infty }({{\mathbb {R}}}),
\end{aligned}\]	
	and
\[\begin{aligned} \int_{{\mathbb {R}}}{{\log }^+}\left( \frac{\left| {\hat{f}}(y)e^{{\frac{{\pi
				}^2}{\alpha }y^2}}\right| }{\rho }\right) \ {\hbox {d}}y<\infty , \end{aligned}
			\]
	for some positive constants $\alpha $ and $\rho $, where
\[\begin{aligned} \log ^+(x)={\left\{ \begin{array}{ll} \log (x), &{} \text {if }x > 1.\\ 0 &{} \text
	{otherwise}. \end{array}\right. } \end{aligned}\]
	
and $\hat{f}$ stands for the classical Fourier transform of $f$,\\
	then $f$ is a constant multiple of the Gaussian $e^{{-\alpha x}^2}$.
	
	Miyachi's theorem has been extended in several different directions in recent years, including extensions to Dunkl transform \cite{CDKM11}, Clifford--Fourier transform\cite{KJ17}, and much more generally, to nilpotent lie groups \cite{BT10} and Heisenberg motion groups \cite{BT16}.
	
The quaternion Fourier transform (QFT) is a non-trivial extension of the real and complex classical Fourier transform to the algebra of the quaternions.
Since the quaternion multiplication is non-commutative, there are three types of the QFT depending on which side multiplication of the kernel is done, that is the so-called left-sided, right-sided and the two-sided QFT.\\
	The QFT was introduced at first by Ell \cite{EL93} for the analysis linear time-invariant partial differential systems and then applied in color image processing.\\
	Later, Bülow \cite{BU99} investigated the important properties of the two-sided QFT for real signals and applied it to signal and image processing.
    Furthermore, several uncertainty principles have been
	formulated for the quaternion Fourier
	transform.
	
    In \cite{MA08}, the authors generalized a component-wise UP
	for the right-sided QFT. The directional UP related to the two-sided QFT was proposed
	in \cite{HI10}. Recently, in \cite{CKL15}, the authors established logarithmic UP associated with the QFT. Meanwhile, Mawardi \cite{MA16} obtained the connection between	the QFT and quantum mechanics and then established the
	modified UP (full UP) for the two-sided QFT.
	
	Our contribution to these developments is that we propose a new UP for the QFT, namely Miyachi's UP.\\
		So far, no such uncertainty principle for the QFT (one-sided or two-sided) had been established. In our previous works, other UPs: Heisenberg, Hardy\cite{EF17}, and Beurling\cite{EF19}, have been extended for the two-sided QFT. Also, we derived in \cite{EF19} the UPs of Cowling-Price and Hardy using the extension of Beurling theorem in a quaternion framework. Here, we will obtain, in a different way, by the main result of Miyachi, the same UPs of Cowling-Price and Hardy in QFT domain.
The techniques used here are also applicable for the left-sided and the right-sided QFT as well.\\
	Our paper is organized as follows. In Sect. 2, we review some basic notions and notations related to the quaternion algebra. In Sect. 3, we recall the definition and some results for the quaternion Fourier transform useful in the sequel. In Sect. 4, we prove Miyachi’s theorem  for the quaternion Fourier transform, and provide an extension of certain UPs to the quaternion Fourier
		transform domain. In Sect. 5, we conclude the paper.

	\section{The Algebra of Quaternions }
	In order to extend complex numbers to a four-dimensional algebra, the Irish W. R. Hamilton invented in 1843 the quaternion algebra $\Q$.\\
	Any quaternion $q\in \Q$\ can be expressed by
	$$q=q_0+\e q_1+\f q_2+\g q_3;\ q_0, q_1, q_2, q_3 \in \R,$$
	where  $\e, \f, \g$  satisfy Hamilton's rules
	
	\bee
	\e^2 = \f^2 = \g^2=-1,\ \e\f = -\f\e = \g,\eee
		 $$\f\g = -\g\f = \e;\ \g\e = -\e\g = \f.$$
	Quaternions are isomorphic to the Clifford algebra ${Cl}_{(0,2)}$ of ${\R}^{(0,2)}$:
\be \label{2.1}
\Q \cong Cl_{(0,2)}.       
\ee

	We define the conjugation of $q \in \Q$\ by:
	$$\overline{q}\myeq q_0-\e q_1-\f q_2-\g q_3$$
and its modulus $|q|_Q$ is defined by
	\[{|q|}_Q\myeq\sqrt{q\overline{q}}=\sqrt{q^2_0+ q^2_1+q^2_2+q^2_3}.\] 
	
Particularly, when $q=q_0$ is a real number, the module ${|q|}_Q$ reduces to the ordinary Euclidean module $\left|q\right|=\sqrt{q^2_0}$.
	Also, we observe that for $\x\in {\R}^2,\ {\left|\x\right|}_Q=\left|\x\right|,$ where $\left|.\right|\ $is the Euclidean norm $\left|(x_1,x_2)\right|^2=x^2_1+x^2_2$\

Moreover, for arbitrary $p,q\in \Q$ the following identity holds
	
	\[{|pq|}_Q={|p|}_Q{|q|}_Q.\] 
	
	Clearly, the inverse of   $0\ne q \in \Q$ is defined by~:
	
	\[q^{-1}=\frac{\overline{q}}{|q|^2_Q}.\] 
which shows that $\Q$ is a normed division algebra.	

	Due to \eqref{2.1}, we recall the following properties:
	
	if   $\x$ is a vector in ${Cl}_{(0,2)}$, then 
	\be \label{2.2}
	{\left|\x\right|}^2=-\x^2. 
	\ee
	Let (,) be the inner product on  ${\R}^{(0,2)}$; for \ref{2.3}rs\ $\x$\ and $\y$\ we have 
	
	\be\label{2.3}
	(\x,\y)~\myeq \sum^2_{l=1}{x_l}y_l=-\frac{1}{2}(\x\y+\y\x).
	\ee
	In this paper, we will study the  quaternion-valued signal $f:\R^2\to \Q $, which can be written in this form $$f=f_0+\e f_1+\f f_2+\g f_3,$$ with $f_m~: {\R^2}\to \ {\R}\ for\ m=0,1,2,3.$\\
	We introduce the Banach spaces $L^p\left({\R}^2,\Q\right)$,  $1\le p \le \infty $, 
	$$L^p\left({\R }^2,\Q\right) \myeq\ \{f|f:\R^2 \rightarrow \Q, {\left|f\right|}_{p,Q}\myeq ({\int_{{\R }^2}{{\left|f(\x)\right|}^p_Q \ d\x}})^{\frac{1}{p}}< \infty \}, 1\le p  <\infty,$$
	 $$L^\infty \left({\R }^2,\Q\right) \myeq\ \{f|f:\R^2 \rightarrow \Q, {\left|f\right|}_{\infty,Q}\myeq  ess \ {sup}_{\x\in \R^2}{\left|f(\x)\right|}_Q< \infty \}.$$
	 where $d\x=dx_1dx_2,$ refers to the usual Lebesgue measure in $\R^2.$
	
	If $f\ {\in L}^{\infty } \left({\R}^2,\Q\right)$  is continuous, then 
		\[{\left|f\right|}_{\infty ,Q}=  {sup}_{{\x\in {\R}}^2}{{\left|f(\x)\right|}_Q}\] 
	Furthermore, we define naturally the two following Banach spaces 
	$$L^1(\R^2,\Q)\cap L^\infty(\R^2,\Q)= \{f|f \in L^1(\R^2,\Q)\  \textnormal{and} \ f\in L^\infty (\R^2,\Q) \},$$
		$$L^1(\R^2,\Q)+ L^\infty(\R^2,\Q)= \{f=f_1+f_2, f_1\in L^1(\R^2,\Q),\ f_2\in L^\infty (\R^2,\Q) \}.$$
	
	We denote by ${\mathcal S}(\R^2,\mathbb H)$ the quaternion Schwartz test function space, i.e., the set $C^{\infty }$ of smooth functions $f$, from ${{ \R}}^2$ to $\mathbb H$, given by 
		$$\mathcal S(\R^2,\Q)\myeq\{f \in C^{\infty}(\R^2,\Q):{sup}_{\x\in \R^2, \ {|\alpha|\le n}}
	{({\left(1+\left|\x\right|^m\right)} \partial^{\mathbf{\alpha} } {\left|f(\x)\right|}_Q)}<\infty,\ m, n\in {\mathbb N}\},$$
	where $\partial^{\mathbf{\alpha} }\myeq \frac{\partial^{|\alpha|}}{\partial^{\alpha_1}_{x_1}\partial^{\alpha_2}_{x_2}}$, $|\alpha|=\alpha_1+\alpha_2$ for a multi-index $\alpha= ({\alpha }_1,{\alpha }_2)\in \mathbb N^2$.	
	\section{Quaternion Fourier Transform}
In this section, we review the definition and some properties of the two-sided QFT.
	\begin{Def}
		Let$\ f$ in    $L^1\left({\R}^2,\Q\right)$. Then, the two-sided quaternion Fourier transform of the function  $f$ is given by
		
		\be \label{3.1}
		\F\{f\}(\xii )\myeq \int_{{\R}^2}{e^{-\e{2\pi \xi }_1x_1}} f(\x) e^{-\f{2\pi \xi }_2x_2}d\x,\ee
		
		where ~$\xii,\x\in {\R}^2.$
	\end{Def}
	\begin{Lemma}\label{inver}{Inverse QFT} \cite[Thm. 2.5]{BU99}
		
		If  $f\in L^2\left(\R^2, \Q\right), and \ \F\{f\}\in L^1\left(\R^2,\Q\right)$, then the two-sided QFT is an invertible transform and its inverse is given
		by
		
		$$f(\x)=\int_{\R^2}{e^{\e2\pi {\xi }_1x_1}} \F\{f\}(\xii ) e^{\f{2\pi \xi }_2x_2}d\xii,\ \  d\xii=d\xi_1d\xi_2. $$
	\end{Lemma}

	\begin{Lemma}\label{sca}{Scaling property} 
		
		Let $\alpha $ be a positive scalar constant; then, the two-sided QFT of \ ${f}_{\alpha }\left({ \x}\right)={ f}(\alpha { \x})$  becomes
\be
\F\left\{f_\alpha\right\}(\xii) ={\left(\frac{1}{\alpha }\right)}^{2}\F\left\{f\right\}(\frac{1}{\alpha }\xii).
\ee
	\end{Lemma}
	
	Proof.  Equation \eqref{3.1} gives
	\bee
	\F\{{{ f}}_{\alpha }\}(\xii )=\int_{{\R}^2}{e^{-\e{2\pi \xi }_1x_1}}\ f(\alpha \x)\ e^{-\f{2\pi \xi }_2x_2}d\x.
	\eee
	We substitute $\y$ for $\alpha \x$ and get 
	\bee
		\F\{{{ f}}_{\alpha }\}(\xii )&=&{\left(\frac1{\alpha }\right)}^2\int_{{\R}^2}{e^{-\e2\pi (\frac{1}{\alpha }{\xi }_1)y_1}} f(\y) e^{-\f2\pi (\frac1{\alpha }{\xi }_2)y_2}d\y\\
		&=&{\left(\frac1{\alpha }\right)}^2\F\left\{f\right\}(\frac{1}{\alpha }\xii).
	\eee
The next lemma states that the QFT of a Gaussian quaternion function is another quaternion Gaussian quaternion function.
	\begin{Lemma}\label{Gauss_QFT}{QFT of a Gaussian quaternion function.}\\
		Consider a two-dimensional Gaussian quaternion function $f$ given by
		$$f(\x)=q e^{-(\alpha_1x_1^2+\alpha_2x_2^2)},$$
		where $q=q_0+iq_1+jq_2+kq_3$\ is a constant quaternion, and $\alpha_1, \alpha_2$\ are positive real constants.\\
		Then
		\be
		\F\{f\}(\xii)=q \frac{\pi}{\sqrt{\alpha_1\alpha_2}}e^{-\frac{\pi^2}{\alpha_1}\xi_1^2-\frac{\pi^2}{\alpha_2}\xi_2^2}.
		\ee
		where ~$\x,\xii \in {\R}^2.$
\end{Lemma}
Proof. Let $g$ be defined by $g(\x) \myeq e^{-(\alpha_1x_1^2+\alpha_2x_2^2)},$\ we have 
$$\F\{f\}(\xii)=q_0\F\{g\}(\xii)+\e  q_1\F\{g\}(\xii)+q_2\F\{g\}(\xii) \f+q_3\e\F\{g\}(\xii)\f,$$
where we used the $\R$-linearity of the QFT and the properties\\  
$\F\{\e h\} = \e\F\{h\}, \F\{\f h\} = \F\{h\}\f,$ and $\F\{\g h\} = \e\F\{h\}\f $\  for $h$ real-valued function.\\
On the other hand, we have
\bee
\F\{g\}(\xii)&=&\int_{{\R}^2}{e^{-\e{2\pi \xi }_1x_1}} g(\x) e^{-\f{2\pi \xi }_2x_2}d\x\\
&=& \int_{\R} e^{-\alpha_1x_1^2} {e^{-\e{2\pi \xi }_1x_1}} dx_1 \int_{\R}e^{-\alpha_2x_2^2} e^{-\f{2\pi \xi }_2x_2}dx_2\\
&=& \sqrt{\frac{\pi}{\alpha_1}}e^{-\frac{\pi^2}{\alpha_1}\xi_1^2}\sqrt{\frac{\pi}{\alpha_2}}e^{-\frac{\pi^2}{\alpha_2}\xi_2^2}\\
&=&\frac{\pi}{\sqrt{\alpha_1\alpha_2}}e^{-\frac{\pi^2}{\alpha_1}\xi_1^2-\frac{\pi^2}{\alpha_2}\xi_2^2}.
	\eee
This completes the proof of Lemma \ref{Gauss_QFT}.\qed

	\begin{Lemma}\label{Beur}\cite[Lemma 3.11]{EF19}
		
		Let $f:\R^2\to \mathbb{R} $  be of the form $$f\left(\x\right)=P(\x)e^{-\pi \alpha {\left|\x\right|}^2},$$ 
		where  $P$ is a polynomial  and $\alpha >0,$
		
		Then      \[\F\{f\}(\xii {\rm )}=\ Q(\xii )\ e^{-\frac{\pi }{\alpha }{\left|\xii \right|}^2},\]
		
	where $Q$ is a quaternion polynomial with $deg P = deg Q.$
	\end{Lemma}
	\section{Miyachi's theorem }
	In this section, we prove Miyachi's theorem for the quaternion Fourier transform. For this, we need the following technical lemma of the complex analysis.

	\begin{Lemma}\label{texnic}\cite[Lemma 1]{CDKM11}	\\
		Let $h$  be an entire function on ${\mathbb C}^2$ such that 
		\[{\left|h\left(\bm{z}\right)\right|}\le {Ae}^{B{\left|Re(\bm{z})\right|}^2}\ for\  all\  \bm{z} \in {\mathbb C}^2\]
		and 
		\[\int_{{\R}^2}{{\log}^+}\left({\left|h\left(\y\right)\right|}\right)d\y<\infty,\] 
		for some positive constants $A$ and $B$.\\
		 Then $h$ is a constant function.

	\end{Lemma}

	\begin{Th}\label{main_result}{(Miyachi’s Theorem)}.\\
		Let $\alpha $, $\beta >0$. Suppose that $f$ is a measurable function such that 
		
		\be \label{4.1} e^{{\alpha \left|\x\right|}^2}f\in L^1({\R}^2,\Q)+ L^{\infty }({\R}^2,\Q),\ee
		
		and   \be \label{4.2}\int_{{\R}^2}{{\log}^+}(\frac{{\left|\F\left\{f\right\}(\y)e^{{\beta \left|\y\right|}^2}\right|}_Q}{\rho })\ d\y < \infty, \ee  
		
		for some $\rho ,\ 0< \rho <+\infty $
		
		Then, three cases can occur:
		
		\begin{enumerate}[label=(\roman*)]
			\item  If $\alpha \beta >{\pi }^2,\ $ then $\ f=0$ almost everywhere.
			
			\item   If $\ \alpha \beta ={\pi }^2,\ $then $f\left(\x\right)= C e^{{-\alpha \left|\x\right|}^2},$  where $C$ is a constant quaternion.
			
			\item  $If\ \alpha \beta <{\pi }^2,\ $then there  exist infinitely many  functions satisfying \eqref{4.1} and \eqref{4.2}.
		\end{enumerate}
	\end{Th}
	
	Proof. \begin{enumerate}[label=(\roman*)]
		\item We first prove the result for the case $\alpha \beta ={\pi }^2.$
	
	By scaling, we can assume that $\alpha =\beta =\pi.$
	
	 Indeed, let $g\left(\x\right)=f(\sqrt{\frac{\pi }{\alpha }}\ \x)$; then, by Lemma \ref{sca} we obtain\\
	\be \label{revision}
	\F\left\{g\right\}\left(\bm{t}\right)=\frac{\alpha }{\pi }\ \F\left\{f\right\}(\sqrt{\frac{\alpha }{\pi  }}\bm{t}), \ \ \ \bm{t}\in \R^2.
	\ee
	In addition, we get by \eqref{4.1}
	\bee 
	e^{{\pi \left|\x\right|}^2}g(\x)=e^{{\pi \left|\x\right|}^2}f(\sqrt{\frac{\pi }{\alpha }}\ \x) \in L^1({\R}^2,\Q)+ L^{\infty }({\R}^2,\Q).
	\eee
	Also, we have
	\bee \int_{\R^2}{{\log}^+}\left(\frac{{\left|\F\left\{f\right\}\left(\y\right)e^{{\beta \left|\y\right|}^2}\right|}_Q}{\rho }\right)d\y &=&\frac{\alpha }{\pi }\int_{{\R}^2}{{\log}^+}\left(\frac{{\left|\F\left\{f\right\}\left(\sqrt{\frac{\alpha }{\pi }}\bm{t}\right)e^{{\frac{\alpha \beta }{\pi }\left|\bm{t}\right|}^2}\right|}_Q}{\rho }\right)d\bm{t}\\
	&\overset{\makebox{\mbox{\normalfont\sffamily \eqref{revision}}}}{=}& \frac{\alpha }{\pi }\int_{{\R}^2}{{\log}^+}\left(\frac{{\left|\F\left\{g\right\}\left(\bm{t}\right)e^{{\pi \left|\bm{t}\right|}^2}\right|}_Q}{\rho ^{'}} \right)d\bm{t} < \infty,
	\eee
where ${\rho }^{'}=\frac{\alpha }{\pi } \rho.  $
	
	If  the result is shown for $\alpha =\beta =\pi ,$
	
	then  $g\left(\x\right)= C e^{{-\pi \left|\x\right|}^2}$, and thus, $f\left( \x\right)=\ g\left(\sqrt{\frac{\alpha }{\pi }}\x\right)= C \ e^{{-\alpha \left|\x\right|}^2}.$
	
	Now, we assume that $\alpha =\beta =\pi :$
	
	Applying the same method as in\ \cite[Thm. 5.3]{EF17}, by complexifying the variable
	$$ \bm{z}=\bm{a}+i_{{\mathbb C}}\bm{b},$$  
	where $\bm{a} =(a_1,a_2),\bm{b}=(b_1,b_2)\in \R^2$,\ and we note by $i_{{\mathbb C}}$\ the complex number which satisfies  $i^2_{\mathbb C}= -1.$
	
	We have $${\left|\bm{z}\right|}^2_Q={\left|\bm{a}\right|}^2+{\left|\bm{b}\right|}^2={\left|\bm{z}\right|}^2,$$
	\ where $\left|.\right|\ $is the Euclidean norm in ${\mathbb C}^2.$
	
	Let \bee
	w\left({ \x}\right)&=&{{ e}}^{\pi ({\left|{{ a}}_1\right|^2{ +}\left|{{ b}}_1\right|^2)}}{{ e}}^{\pi ({\left|{{ a}}_2\right|^2{ +}\left|{{ b}}_2\right|^2)}}{{ e}}^{{ -}\pi ({\left|x_1\right|-(\left|{{ a}}_1\right|{ +}\left|{{ b}}_1\right|))}^2}{{ e}}^{{ -}\pi ({\left|x_2\right|-(\left|{{ a}}_2\right|{ +}\left|{{ b}}_2\right|))}^2}\\
    &=&e^{\pi {\left|{ z}\right|}^2}{{ e}}^{{ -}\pi ({\left|x_1\right|-(\left|{{ a}}_1\right|{ +}\left|{{ b}}_1\right|))}^2}{{ e}}^{{ -}\pi ({\left|x_2\right|-(\left|{{ a}}_2\right|{ +}\left|{{ b}}_2\right|))}^2}.
    \eee
	Clearly, $w$\ belongs to $L^1(\R^2,\Q)\cap L^{\infty }(\R^2,\Q).$
	
	By assumption, $e^{{\pi \left|\x\right|}^2}f$ belongs to $L^1$(${{ \R}}^2{ ,\Q)+}L^{\infty }$(${{ \R}}^2{ ,\Q)}.$
	
	As $$\F\left\{f\right\}\left(\bm{z}\right)=\int_{{{ \R}}^2}{e^{-2\pi \e x_1{ (}a_1{ +}i_{\mathbb C} b_1)}f\left(\x\right)e^{-2\pi \f{ \x}_2(a_2+i_{\mathbb C} b_2)}d\x},$$
	we have 
	\bee
	\left| \F\left\{f\right\}\left(\bm{z}\right)  \right|_Q     &\le& \int_{{{ \R}}^2}{{\left|f\left(\x\right)\right|}_Q}e^{2\pi (\left|x_1a_1\right|+\left|x_1b_1\right|{ +}\left|x_2a_2\right|+\left|x_2b_2\right|)}d\x\\
    &=&\int_{\R^2}{{\left|e^{{\pi \left|\x\right|}^2}f\left({ \x}\right)\right|}_Q\ }{ w}\left({ \x}\right)d\x.
    \eee 
	
	Hence, $\F\left\{{ f}\right\}\left(\bm{z}\right)$  is well defined,
	and is an entire function on ${\mathbb C}^2.$
	
	Furthermore, by \eqref{4.1} there exists $u \in L^1({{ \R}}^2{ ,\Q)}$ and $v \in L^{\infty }({{ \R}}^2{ ,\Q)}$  such that
	
	\[e^{{\pi \left|\x\right|}^2}f\left({ \x}\right)= u\left({ \x}\right)+ v\left({ \x}\right),\] 
	Using the triangle inequality and the linearity of the integral we get\\
	\[{\left|\F\left\{f\right\}\left(\bm{z}\right)\right|}_Q\le \int_{\R^2}{{\left|u\left({ \x}\right)\right|}_Q }{ w}\left({ \x}\right) d{ \x}+\int_{\R^2}{{\left|v\left({ \x}\right)\right|}_Q\ }{ w}\left({ \x}\right) d{ \x}.\] 
	
	Then according to the H\^older's inequality, we have
	
	\[{\left|\F\left\{f\right\}\left(z\right)\right|}_Q\le {\left|u\right|}_{1,Q} {\left|w\right|}_{\infty ,Q} + {\left|v\right|}_{\infty ,Q}{\left|w\right|}_{1,Q},\] 

	Since \[\int_{\R}{e^{-\pi ({\left|t\right|+m)}^2}}dt\le 2,\ \ \textnormal{where}\ m\in {\R},\]
	we obtain
	$${\left|w\right|}_{1,Q}\le 4e^{\pi {\left|\bm{z}\right|}^2},\ \textnormal{and}\ {\left|w\right|}_{\infty ,Q}\le  e^{\pi {\left|\bm{z}\right|}^2}.$$
	
	Then \bee {\left|\F\left\{f\right\}\left(\bm{z}\right)\right|}_Q &\le& {e^{\pi {\left|\bm{z}\right|}^2}\left|u\right|}_{1,Q}+4\ e^{\pi {\left|\bm{z}\right|}^2}{\left|v\right|}_{\infty ,Q}\\
	&\le& K\ e^{\pi {\left|\bm{z}\right|}^2},
	\eee
where $K$ is a positive constant independent of $z.$
	
Now, let  $$h(\bm{z})= e^{-\pi z^2}|\F\left\{f\right\}\left(\bm{z}\right)|_Q,\	\textnormal{for} \ \bm{z}\in {\mathbb C}^2,$$ then, $h$ is an entire function.
	
	By \eqref{2.2} and \eqref{2.3}, we have 
	
	\[\bm{z}^2=\left(\bm{a}+i_{\mathbb C} \bm{b}\right)\left(\bm{a}+i_{\mathbb C}\bm{b}\right)=-{\left|\bm{a}\right|}^2+{\left|\bm{b}\right|}^2-2i_{\mathbb C} (\bm{a},\bm{b}),\] 
	then
$${\left|e^{-\pi \bm{z}^2}\right|}_Q\le e^{\pi {\left|\bm{a}\right|}^2}e^{-\pi {\left|\bm{b}\right|}^2},$$
	where we used ${\left|e^{2\pi i_{\mathbb C}(\bm{a},\bm{b})}\right|}_Q=1$.\\
	As a result\\
$$\left|h\left(z\right)\right|\le {Ke}^{\pi {\left|a\right|}^2}e^{-\pi {\left|b\right|}^2}  e^{\pi {\left|a\right|}^2}e^{\pi {\left|b\right|}^2},$$
	   \be \label{3.3}
	   =Ke^{2\pi {\left|a\right|}^2}.
	   \ee
	On the other hand, 
	by   \eqref{2.2} and by assumption
		\be \label{4.4}
	\int_{{\R}^2}{{\log}^+}\left(\frac{\left|h\left(\y\right)\right|}{\rho }\right)d\y=\int_{{\R}^2}{{\log}^+}\left(\frac{{\left|e^{\pi {\left|\y\right|}^2}\ \F\left\{f\right\}\left(\y\right)\right|}_Q}{\rho }\right)d\y<\infty .
	\ee 
		Then by \eqref{3.3},\eqref{4.4} and by applying Lemma \ref{texnic} to the function $h\left(\y\right){ /}\rho $, we deduce that $h\left(\y\right){ =const.}$
	
	i.e  $$|\frac{\F\left\{f\right\}\left(\y\right)}{e^{-\pi {\left|\y\right|}^2}}|_Q=\ const.$$
	Therefore $$\F\left\{f\right\}\left(\y\right)= C e^{-\pi {\left|\y\right|}^2},$$
	where $C$ is a constant quaternion.\\	
	Then by Lemmas \ref{inver}\ and \ref{Gauss_QFT}, we have 
	\[f\left(\x\right) = C\ e^{-\pi {\left|\x\right|}^2}.\] 
	\item If $\alpha \beta >{\pi }^2.$            \\
	Let  $g\left(\x\right)=f(\sqrt{\frac{\pi }{\alpha }}\ \x)$, a simple calculation shows that  
	\bee
	\int_{\R^2}{{\log}^+}\left(\frac{{\left|\F\left\{g\right\}\left(t\right)e^{{\pi \left|t\right|}^2}\right|}_Q}{\rho^{'} }\right)dt &<& \frac{\pi }{\alpha }\int_{{\R}^2}{{\log}^+}\left(\frac{{\left|\F\left\{f\right\}\left(\y\right)e^{{\beta \left|\y\right|}^2}\right|}_Q}{\rho}\right)d\y\\
	&<& \infty, \ \ \ \ \ \ (\textnormal{by}\ \eqref{4.2})\eee
	where ${\rho }^{'}=\frac{\alpha }{\pi }\ \rho.$\\
	Then, according to the first case 
	$$g\left(\x\right) = C\ e^{-\pi {\left|\x\right|}^2},$$ 
	where $C$ is a constant quaternion.\\
	Consequently
	 $${ f}\left({ \x}\right)= C\ e^{-\alpha {\left|\x\right|}^2}.$$
	
	Hence, by Lemma \ref{Gauss_QFT} we get\ $$\F\left\{f\right\}\left(\y\right)=C e^{-\frac{{\pi }^2}{\alpha }{\left|\y\right|}^2}.$$
	Refering to \eqref{4.2}, $C$\  must be zero.
	
	\item For the final case, $\alpha \beta <{\pi }^2$
	
	Let ${f\left(\x\right)=\varphi }_{k,l}(\x)\ e^{-\pi \gamma {\left|\x\right|}^2}$   with  $\frac{\alpha }{\pi }$$<$$\gamma $$<$$\frac{\pi }{\beta },$
	
	where ${{\{\varphi }_{k,l}\}}_{k,l\in {\mathbb N}}$ is a basis of\ \  ${\mathcal S}({\R}^2,\mathbb H)$,  which is defined by 
	$${\varphi }_{k,l}\left( x_1,x_2\right)\myeq {\varphi }_k\left(x_1\right){\varphi }_l\left(x_2\right),$$ 
	for $\left( x_1,x_2\right)\in \R^2,$
	
	and 
$${\varphi }_k\left(x\right)=\frac{{(-1)}^k}{k!}e^{\pi x^2}\frac{d^k}{{dx}^k}(e^{-{ 2}\pi \x^2}),  x\in \R.$$

	It is important to see that${{\{\varphi }_{k,l}\}}_{k,l\in \mathbb N}$ is a basis of ${\mathcal S}(\R^2,\mathbb H)$  (see \cite{DB12, EF17}).\\
	We have \bee
	e^{\alpha {\left|\x\right|}^2}f &=& e^{\alpha {\left|\x\right|}^2}{\varphi }_{k,l}\left(\x\right)e^{-\pi \gamma {\left|\x\right|}^2}\\
	&=&{\varphi }_{k,l}\left(\x\right)e^{(\alpha -\pi \gamma ){\left|\x\right|}^2}\in L^1(\R^2{ ,\Q) + }L^{\infty }(\R^2,\Q).\eee
   Lemma \ref{Beur} implies 
	$$\F\left\{f\right\}\left(\y\right)=Q(\y)\ e^{-\frac{\pi }{\gamma }{\left|\y\right|}^2},$$
	where $Q$\ is a quaternion polynomial.
	
	Then, since $\beta <\frac{\pi }{\alpha }$, 
	
	\bee
	\int_{\R^2}{{\log}^+}\left(\frac{{\left|\F\left\{f\right\}\left(\y\right)e^{{\beta \left|\y\right|}^2}\right|}_Q}{\rho }\right)d\y   &=& \int_{\R^2}{{\log}^+}\left(\frac{{\left|Q(\y)\ e^{{(\beta -\frac{\pi }{\gamma })\left|\y\right|}^2}\right|}_Q}{\rho }\right)d\y\\
	&<&\infty.
	\eee
\end{enumerate} 
	This completes the proof of Theorem \ref{main_result}.\qed\\
In the following,  we illustrate the effectiveness of Theorem \ref{main_result}, by giving an example, and derive two generalizations of uncertainty principle associated with the QFT.
 \begin{Ex}
Consider $\alpha, \beta$ two positive numbers with $\alpha \beta =\pi^2,$ and the quaternion Gaussian function 
$$f(\x)=qe^{-\alpha |\x|^2},$$
where $q=q_0+\e q_1+\f q_2+\g q_3$\ is a constant quaternion.\\
Obviously, we have $$e^{\alpha |\x|^2}f=q \in L^\infty(\R^2,\Q)\subset L^1(\R^2,\Q)+L^\infty(\R^2,\Q).$$
Thus, $f$\ satisfies condition \eqref{4.1}.\\
By lemma \ref{Gauss_QFT}, we get $\F\{f\}(\y)=q\frac{\pi}{\alpha}e^{-\beta|\y|^2}.$\\
Then
$$\int_{{\R }^2}\log^{+}(\frac{|\F\{f\}(\y)|_Q e^{\beta|\y|^2}}{\rho})d\y=\int_{{\R }^2}\log^{+}(\frac{|q|_Q\frac{\pi}{\alpha}}{\rho})d\y<\infty,$$
whenever $\rho >|q|_Q\ \frac{\pi}{\alpha}.$
\end{Ex}

\begin{Cor}{Hardy's uncertainty principle for the QFT.}\\
Let $\alpha$ and $\beta$ be both positive constants. Suppose  $f$\ be in $L^1(\R^2,\Q)$
\ with
\begin{enumerate}[label=(\roman*)]
	\item $|f(\x)|^2 < C e^{-\alpha |\x|^2}.$
    \item $|\F \{f\} (\y)|^2 < C^{'} e^{-\beta |\y|^2}.$ 
\end{enumerate}
	for some constants $C>0$ and $C'>0$. Then, three cases can occur :
	\begin{enumerate}
		\item If $\alpha \beta > \pi^2,$ \ $f=0$\ almost everywhere on $\R^2.$ 
		\item If $\alpha \beta = \pi^2,$ then $f$ is a constant quaternion multiple of $e^{-\alpha |\x|^2}$.
		\item  If $\alpha \beta < \pi^2,$ there are infinitely many linearly independent functions satisfying both conditions $(i)$ and $(ii).$
	\end{enumerate}
\end{Cor}
Proof. Immediately using the decay condition $(i)$ one has $f e^{\alpha |\x|^2}\in L^\infty(\R^2,\Q).$\ Hence $f$ verifies condition \eqref{4.1} of Theorem \ref{main_result}.\\
Moreover, for $\rho > 0$\ we have
\[\frac{|\F\{f\}(\y)|_Q e^{\beta|\y|^2}}{\rho}\le\frac{C^{'}}{\rho}.\]
Thus
$$\log^{+}(\frac{|\F\{f\}(\y)|_Q e^{\beta|\y|^2}}{\rho}) \le \underbrace{\log^+(\frac{C^{'}}{\rho})}_{=0},$$
whenver $\rho > C^{'}.$\\
So condition \eqref{4.2} of Theorem \ref{main_result} is verified.
Then, direct application of Theorem \ref{main_result} enables us to achieve the proof.
\begin{Cor}{Cowling-Price's uncertainty principle for the QFT.}\\
   Let	$\alpha$ and $\beta$ be positive real numbers, $1\le p,q\le \infty$ such that $min(p,q)$\ is finite, and let $f$ are a square integrable quaternion-valued function satisfying the following decay conditions  Suppose  $f$\ be in $L^1(\R^2,\Q)$
	\begin{enumerate}[label=(\roman*)]
		\item $\int_{\R^2} (|f(\x)|_Qe^{\alpha |\x|^2})^pd\x < \infty.$
		\item $\int_{\R^2} (|\F\{f\}(\y)|_Qe^{\beta |\y|^2})^q d\y < \infty.$ 
	\end{enumerate}
	Then the three following conclusions hold:
	\begin{enumerate}
		\item $f=0$\ almost everywhere whenever $\alpha \beta > \pi^2.$ 
		\item If $\alpha \beta = \pi^2,$ then $f$ is a constant quaternion multiple of $e^{-\alpha |\x|^2}$.
		\item  If $\alpha \beta < \pi^2,$ there are infinitely many linearly independent functions satisfying both conditions $(i)$ and $(ii).$
	\end{enumerate}
\end{Cor}
Proof. According to (i), we get $f e^{\alpha |\x|^2}\in L^p(\R^2,\Q)$, and using the fact that $L^p\subset L^1+L^\infty$ we obtain that $f$ fulfills condition \eqref{4.1} of Theorem \ref{main_result}.\\ 
Furthermore, based on the inequality
\bee \label{log_prop}
\log^{+}(x)\le x  \ \  \textnormal{for}\ x \in \R_{+},
\eee
 we can easily see that $f$ satisfies the second condition \eqref{4.2}.\\
  Therefore, by applying Theorem \ref{main_result} we conclude the proof.
	\section{Conclusions}
	In this paper, based on some obtained results of the two-sided QFT and one technical  lemma of the complex analysis, a generalization
	of Miyachi's uncertainty principle associated with the QFT was proposed. Consequently, two variants of this UP were provided, namely the theorems of Hardy and Cowling-Price.
	The extension of these qualitative UPs to the quaternionic algebra framework shows that
	a quaternionic 2D signal $f$\ and its QFT cannot both simultaneously decrease very rapidly.\\
The QFT has proved to be very significant tool for applications in color image processing\cite{HKK18}, quantum mechanics, engineering, signal processing, optics, etc. Apart from their importance to pure mathematics, our results are also relevant to applied mathematics and signal processing.


\end{document}